\documentclass[12pt]{article}

\usepackage{amsmath,amsthm}
\usepackage{amssymb}
\usepackage{color}
\usepackage{xspace}
\usepackage{graphicx}
\usepackage[colorlinks=true,
linkcolor=green,
filecolor=brown,
citecolor=green]{hyperref}

\def\v{\vert}

\def\lra{left-to-right maxima\xspace}
\def\lrm{left-to-right maximum\xspace}

\def\j{k}
\def\k{j}

\def\a{\ensuremath{\mathcal A}\xspace}

\def\p{\ensuremath{\mathcal P}\xspace}

\def\t{\ensuremath{\mathcal T}\xspace}

\newcommand{\StirlingPartition}[2]{\genfrac{ \{ }{ \} }{0pt}{}{#1}{#2}}

\newcommand{\seqnum}[1]{\href{http://oeis.org/#1}{\underline{#1}}}

\begin{document}
\newtheorem{theorem}{Theorem}
\newtheorem{defn}[theorem]{Definition}
\newtheorem{lemma}[theorem]{Lemma}
\newtheorem{prop}[theorem]{Proposition}
\newtheorem{cor}[theorem]{Corollary}
\begin{center}
{\Large
The Number of $\overline{3}\overline{1}542$-Avoiding Permutations    \\ 
}

\vspace*{5mm}

DAVID CALLAN  \\
Department of Statistics  \\
\vspace*{-1mm}
University of Wisconsin-Madison  \\
\vspace*{-1mm}
1300 University Ave  \\
\vspace*{-1mm}
Madison, WI \ 53706-1532  \\
{\bf callan@stat.wisc.edu}  \\
\vspace*{5mm}
\end{center}

\begin{abstract}
We confirm a conjecture of Lara Pudwell and show that 
permutations of $[n] $ that avoid the barred pattern 
$\overline{3}\overline{1}542$ are counted by OEIS 
sequence A047970. In fact, we show bijectively that the number of $\overline{3}\overline{1}542$ 
avoiders of length $n$ with $j+k$ left-to-right maxima, of which $\k$ 
initiate a descent in the permutation and $\j$ do not, is 
$\binom{n}{\j}\,\k!\,\StirlingPartition{n-j-k}{\k}$, where 
$\StirlingPartition{n}{\k}$ is the Stirling partition number.

\end{abstract}

\vspace{5mm}

\section{Introduction} \vspace*{-3mm}
A permutation $\pi$ avoids the barred pattern $\overline{3}\overline{1}542$ 
if each instance of a not-necessarily-consecutive 542 pattern in $\pi$ 
is part of a 31542 pattern in $\pi$.
Lara Pudwell \cite{schemes} introduces an automated method to produce ``prefix enumeration 
schemes'' for permutations avoiding a set of one or more barred patterns. Such a scheme, if found, 
translates to a recurrence to count the corresponding avoiders. For single patterns of length 
5 with 2 bars, she notes a success rate of $136/172\ (79.1\%)$. The title pattern is 
among the failures but Pudwell observes \cite{comment1} that the counting sequence 
for this pattern appears to be \seqnum{A047970} in The On-Line Encyclopedia of 
Integer Sequences \cite{oeis}. Here we confirm this conjecture. The definition of 
A047970 is $a(n)=\sum_{i=0}^{n} \big((i + 1)^{n - i} - i^{n - i}\big)$, 
which is equivalent to $a(n)=\sum_{j,k\ge 0}\binom{n}{\j}\,\k!\,
\StirlingPartition{n-j-k}{\k}$, where 
$\StirlingPartition{n}{\k}$ is the Stirling partition number. (This 
equivalence is not entirely obvious and a proof is outlined in the 
Appendix.) Our result is the following.
\begin{theorem}\label{thm1}
  The  number of $\overline{3}\overline{1}542$ 
avoiders of length $n$ with $j+k$ left-to-right maxima, of which $\k$ 
initiate a descent in the permutation and $\j$ do not, is 
$\binom{n}{\j}\,\k!\,\StirlingPartition{n-j-k}{\k}$.
\end{theorem}
Section 2 explains the $\binom{n}{\j}$ factor, and reduces the problem to the case $\j=0$. Section 3 looks at the structure of the avoiders with $\j=0$ and presents a bijection to complete the proof.
To standardize a permutation means to replace its smallest entry by 1, next smallest by 2, and so on.
\vspace*{-5mm}
\section{Reduction} \vspace*{-3mm}
We begin with a reduction of the problem. Deleting the \lra that do not 
initiate a descent and 
standardizing the remaining permutation does not affect the avoider status 
because a \lrm can serve only as the 5 of a 542 pattern or as the 3 
in the 31 of a 31542 pattern, and each 
deleted \lrm is immediately followed (in the original permutation) by another \lrm or ends the path and so ``it won't be missed''. 
Furthermore, the original permutation can be recovered from knowledge 
of the deleted entries and the standardized permutation. These observations
account for the $\binom{n}{\j}$ factor (choose $\j$ entries from $[n]$ 
to serve as the \lra that do not initiate a descent), and it
means that Theorem \ref{thm1} will follow from 
\begin{theorem}\label{thm2}
    The number of $\overline{3}\overline{1}542$-avoiding permutations 
    with $\k$ \lra, all of which initiate a descent, is $\k!\,\StirlingPartition{n-\k}{\k}$.
\end{theorem}
\vspace*{-5mm}
\section{Proof of Theorem 2} \vspace*{-3mm}
To characterize $\overline{3}\overline{1}542$-avoiders in 
which  each \lrm initiates a descent, it is 
convenient to represent a permutation $\pi$ as the usual matrix diagram, as in Figure 1, 
with  a bullet in the $(a,b)$ position (measuring from lower left) if and only 
if $\pi(a)=b$. 

\begin{figure}
\vspace{-.6in}
\begin{center}
\includegraphics[angle=0, scale = 1]{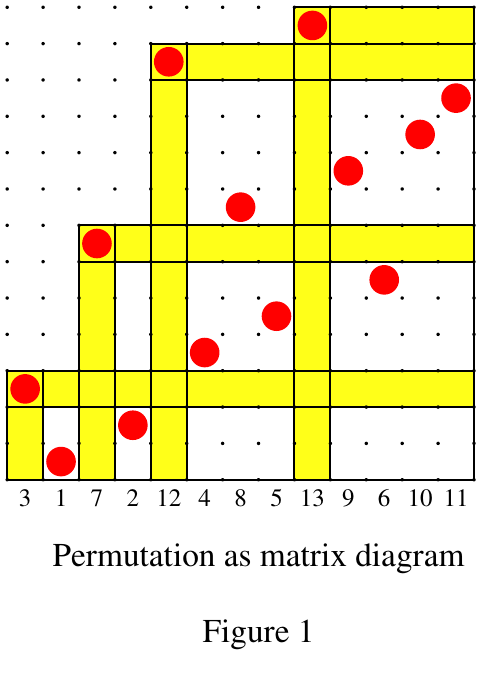}
\end{center}
\end{figure}
The  \lra determine the upper left staircase boundary 
and the yellow lines through the \lra delineate a collection of
$1+2+\ldots+\k$ cells (the white rectangles), where $\k$ is the number of \lra. 
Some cells may be empty (contain no bullets) or even vacuous (contain no area).
The cells in turn split into horizontal strips and also into vertical strips: 
the $i$-th horizontal (resp.\ vertical) strip contains  $\k+1-i$ (resp.\:$i$) 
cells. The condition that  each \lrm initiates a descent means 
that no vertical strip consists entirely of empty cells. 
\begin{prop}
  Suppose each \lrm in a permutation $\pi$ initiates a descent. Then  
  $\pi$ is $\overline{3}\overline{1}542$-avoiding if and only if the 
  bullets in each horizontal strip are rising 
  from left to right. 
\end{prop}

\vspace*{-2mm}
Proof. 
Two falling bullets correspond to an inversion 
$ba$ in the permutation. 
If they lie in the same horizontal strip, let $c$ denote the \lrm associated with 
this strip. Then $cba$ is a 542 with no available 3 let alone 31. 

On the 
other hand, neither the 4 nor 2 of a 542 pattern can be a \lrm and, if 
the rising bullet condition holds, the 4 and 2 lie in \emph{different} 
horizontal strips. The \lrm $m$ associated with the horizontal 
strip containing the 2 serves as the 3. The vertical strip associated 
with $m$ is nonempty (since all vertical strips are nonempty) and any 
bullet in it serves as the 1.  \qed

Let $\t_{\k}=[1]\times [2]\times \ldots \times [\k]$ (which may conveniently be referred to as 
the set of inversion codes for permutations of $[\k]$), 
and let $\p_{n,\k}$ denote the set of partitions of $[n]$ 
into $\k$ blocks in a canonical form: the smallest entry of each block 
is last, the remaining entries increase left to right, and the blocks 
are arranged in order of increasing smallest entry. For example, $1\,/\,
5\,7\,2\,/\,4\,8\,9\,3\,/\,6 \in \p_{9,4}$. Clearly, $\v\,\t_{\k}\,\v = 
\k!$ and  $\v\,\p_{n,\k}\,\v = \StirlingPartition{n}{\k}$.

The \emph{successors} in 
$\pi$ are the entries immediately following the \lra, that is, the 
terminators of the descents which the \lra initiate.
Here is a bijection from the set $\a_{n,\k}$ of 
$\overline{3}\overline{1}542$-avoiding permutations of $[n]$ with $\k$ 
\lra, all of which initiate descents, to $\t_{\k}\times \p_{n-\k,\,\k}$.

Given $\pi \in \a_{n,\k}$, draw its diagram. Record the horizontal strips 
(numbered bottom to top) of 
the successors in $\pi$. In the example of Figure 1, the successors are $(1,2,4,9)$ and 
their horizontal strips are $(1,1,2,3)$. This is the desired 
inversion code. 

\begin{figure}
\vspace{-.6in}
\begin{center}
\includegraphics[angle=0, scale = 1]{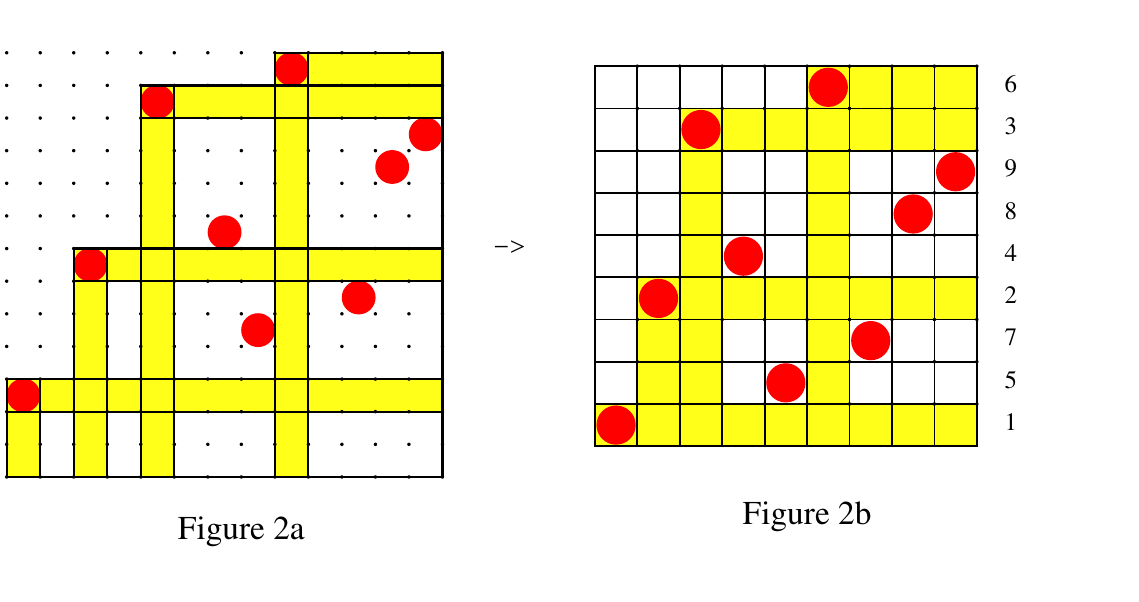}
\end{center}
\end{figure}
Erase the bullets of the successors (Figure 2a), ``prettify'' 
the diagram (Figure 2b), record the 
position of each bullet, bottom to top, among all the  
bullets and place a divider to the  right of each right-to-left minimum. 
This is the desired partition in 
$\p_{n-\k,\,\k}$ and it is already in canonical form. 
The example yields 
$1\,/\,5\,7\,2\,/\,4\,8\,9\,3\,/\,6$. 

The map is reversible because the diagram of Figure 2b can be recovered from the partition, the inversion code determines the cell 
of each missing bullet and, within the cell, the bullet's horizontal position 
is the extreme left (since it lies immediately after a \lrm) and its vertical 
position is determined by the rising bullet condition.
\vspace*{-5mm}

\section{Appendix} \vspace*{-3mm}
\begin{prop}
\begin{equation}
  \sum_{i=0}^{n} \big((i + 1)^{n - i} - i^{n - i}\big)=
 \sum_{j,k\ge 0}\binom{n}{\j}\,\k!\,\StirlingPartition{n-j-k}{\k}.   
    \label{identity}
\end{equation}
\end{prop}
Proof. Let $A_{n}= \sum_{i=0}^{n}  i^{n - i}$. Then  
$\sum_{i=0}^{n}  (i+1)^{n - i}=A_{n+1}$, and the left hand side of 
(\ref{identity}) is $A_{n+1}-A_{n}$.
The binomial theorem says 
$A_{n+1}=\sum_{i=0}^{n}\sum_{\ell=0}^{n-i}\binom{n-i}{\ell}i^{\ell}$, 
leading to 
\[
A_{n+1}-A_{n}=\sum_{\ell=0}^{n}\sum_{i=0}^{n-\ell}\left(\binom{n-i}{\ell}-
\binom{n-i-1}{\ell}\right)i^{\ell} =
\sum_{\ell=0}^{n}\sum_{i=0}^{n-\ell}\binom{n-i-1}{\ell-1}i^{\ell}.
\]

On the other hand, the right side of (\ref{identity}) is
\begin{eqnarray*}
    \sum_{j,k\ge 0}\binom{n}{\j}\,\k!\,
\StirlingPartition{n-j-k}{\k} & \overset{\textrm{(i)}}{=} & 
\sum_{\j=0}^{n}\sum_{\k=0}^{n-\j}\sum_{i=0}^{\k}\binom{n}{\j}(-1)^{\k-i}\binom{\k}{i}i^{n-j-k}  \\
     & \overset{\textrm{(ii)}}{=} & 
     \sum_{\ell=0}^{n}\sum_{i=0}^{\ell}i^{n-\ell}\sum_{\j=0}^{\ell-i} 
     \binom{n}{\j}(-1)^{\ell-i-\j}\binom{\ell-\j}{i}  \\
     & \overset{\textrm{(iii)}}{=} &   \sum_{\ell=0}^{n}\sum_{i=0}^{\ell}i^{n-\ell} 
     \binom{n-i-1}{n-\ell-1} \\
     & \overset{\textrm{(iv)}}{=} & 
     \sum_{\ell=0}^{n}\sum_{i=0}^{n-\ell}i^{\ell}\binom{n-i-1}{\ell-1}, 
\end{eqnarray*}
agreeing with the left side.

\noindent Notes on equalities: 
\vspace*{-1mm}
\begin{enumerate}

\item Expand Stirling partition number.
\vspace*{-2mm}
\item Change summation index from $\k$ to $\ell$  with $\ell=j+k$.
\vspace*{-2mm}
\item  Evaluate inner sum; this is the essential content of identity
(5.25) in \cite{gkp}. The identity (5.25) in \cite{gkp} ostensibly
has 4 parameters, but a change of summation index from $ k$ to
$j=k-n$ in (5.25) leaves only 3 independent parameters, as here.
\vspace*{-2mm}
\item Reverse order of outer summation.

\end{enumerate}

\end{document}